\documentclass{article}
\usepackage{amsmath,amssymb,latexsym,amsthm,amscd}
\font\goth=eusm10

\newcommand\bc{\mathbf C}

\newcommand\E{\mathcal E}

\newcommand\bz{\mathbf{Z}}

\newcommand\bn{\mathbf{N}}

\newcommand\Oc{\hbox{\goth O}}

\newtheorem{theorem}{Theorem}

\newtheorem{corollary}{Corollary}

\newtheorem{lemma}{Lemma}

\numberwithin{proposition}{section}
\numberwithin{definition}{section}
\numberwithin{corollary}{section}
\numberwithin{remark}{section}
\numberwithin{lemma}{section}
\numberwithin{equation}{section}

\numberwithin{theorem}{section}

\numberwithin{question}{section}
\numberwithin{case}{section}
\numberwithin{example}{section}
\numberwithin{conjecture}{section}

\begin{document}
\title{On correspondences of a K3 surface\\ with itself. III}
\author{C.G.Madonna \footnote{Supported by I3P contract} \ and
Viacheslav V.Nikulin \footnote{Supported by EPSRC grant
EP/D061997/1}}
\maketitle

\begin{abstract} Let $X$ be a K3 surface, and $H$ its primitive polarization
of the degree $H^2=2rs$, $r,s\ge 1$. The moduli space of sheaves
over $X$ with the isotropic Mukai vector $(r,H,s)$ is again a K3
surface, $Y$. In \cite{Mad-Nik1}, \cite{Mad-Nik2} and \cite{Nik1}
(in general) we gave necessary and sufficient conditions in terms
of Picard lattice $N(X)$ of $X$ when $Y$ is isomorphic to $X$,
under the additional condition $H\cdot N(X)=\bz$.

Here we show that these conditions imply existence of an
isomorphism between $Y$ and $X$ which is a composition of some 
universal isomorphisms between moduli of sheaves over $X$,
and Tyurin's isomorphsim between moduli of sheaves over $X$ and
$X$ itself. It follows that for a general K3 surface $X$ with $H\cdot N(X)=\bz$ 
and $Y\cong X$, there exists an isomorphism $Y\cong X$ which is a composition 
of the universal and the Tyurin's isomorphisms. 

This generalizes our recent results \cite{Mad-Nik3}
for $r=s=2$ on similar subject.
\end{abstract}

\section{Introduction} \label{introduction}
We consider algebraic K3 surfaces $X$ over the field $\bc$ of
complex numbers. For a Mukai vector $v=(r,c_1,s)$ where $r\in
\bn$, $s\in \bz$ and $c_1\in N(X)$, Picard lattice of $X$, we
denote by $Y=M_X(r,c_1,s)$ the moduli space of stable (with
respect to some ample $H^\prime\in N(X)$) rank $r$ sheaves on $X$
with first Chern classes $c_1$, and Euler characteristic $r+s$.
The general common divisor of the Mukai vector $v$ is
$$
(r,c_1,s)=(r,d,s)
$$
if $c_1=dc_1^\prime$ where $d\in\bn$ and $c_1^\prime$ is primitive
in the Picard lattice $N(X)$ which is a free $\bz$-module of the
rank $\rho(X)=\text{rk } N(X)$. Here $c_1^\prime\in N(X)$ is {\it
primitive} means that $c_1^\prime/n\notin N(X)$ for any natural
$n\ge 2$. A Mukai vector $v$ is called {\it primitive} if its
general common divisor is one.

By results of Mukai \cite{Muk1}, \cite{Muk2},
under suitable conditions on the Chern classes, the moduli space $Y$ is
always deformations equivalent
to a Hilbert scheme of 0-dimensional cycles on $X$ (of the same dimension).
\par\medskip

In \cite{Mad-Nik1} for $r=s=2$, in \cite{Mad-Nik2} for $r=s=c\in
\bn$, and in \cite{Nik1} for arbitrary   $r,s\in \bn$, the
following result had been obtained.

\begin{theorem}
\label{maintheoremint} Let $X$ be a K3 surface with a polarization
$H$ such that $H^2=2rs$, $r,s\ge 1$, the Mukai vector $(r,H,s)$ be
primitive, and $Y=M_X(r,H,s)$ be the moduli of sheaves on $X$ with
the Mukai vector $(r,H,s)$.

Then $Y\cong X$ if at least for one of signs $\pm$ there exists
$h_1 \in N(X)$ such that the elements $H$ and $h_1$ are contained
in a 2-dimensional sublattice $N \subset N(X)$ with $H\cdot N=\bz$
and $h_1$ belongs to either the $a$-series or the $b$-series
described below where $c=(r,s)$, $a=r/c$, $b=s/c$:

$a$-series:
\begin{equation}
h_1^2=\pm 2bc, \ \ H\cdot h_1\equiv 0\mod bc,
\label{aseriesint}
\end{equation}

$b$-series:
\begin{equation}
h_1^2=\pm 2ac, \ \ H\cdot h_1\equiv 0\mod ac.
\label{bseriesint}
\end{equation}

The conditions above are necessary for $H\cdot N(X)=\bz$ and
$Y\cong X$ if $\rho (X)\le 2$ and $X$ is a general K3 surface with
its Picard lattice, I. e., the automorphism group of the
transcendental periods $(T(X),H^{2,0}(X))$ of $X$ is $\pm 1$.
\end{theorem}

In formulation of this Theorem \ref{maintheoremint} in \cite{Nik1} some
additional primitivity conditions on $h_1$ were also required.
We will show in Sect. \ref{section2}, Remark 2.1,
that they are unnecessary.

The sufficient part of the proof of Theorem \ref{maintheorem} in
\cite{Nik1} and similar Theorems in \cite{Mad-Nik1},
\cite{Mad-Nik2} used global Torelli Theorem for K3 surfaces
\cite{PShShaf}. I. e., under conditions of Theorem
\ref{maintheorem}, we proved that the K3 surfaces $X$ and $Y$ have
isomorphic periods. By global Torelli Theorem,  then $X$ and $Y$
are isomorphic.

In Sect. \ref{section3} below we will give a geometric construction of the isomorphism
between $X$ and $Y$ which is similar to
our considerations in \cite{Mad-Nik3}. This is the main result of this paper.

We prove the following result.

\begin{theorem}
\label{maintheorem2int} Let $X$ be a K3 surface with a
polarization $H$ such that $H^2=2rs$, $r,s\ge 1$, the Mukai vector
$(r,H,s)$ be primitive, and $Y=M_X(r,H,s)$ be the moduli of
sheaves on $X$ with the Mukai vector $(r,H,s)$.

Assume that at least for one of signs $\pm$ there exists $h_1 \in
N(X)$ such that the elements $H$ and $h_1$ are contained in a
2-dimensional sublattice $N \subset N(X)$ with $H\cdot N=\bz$ and
$h_1$ belongs to either the $a$-series or the $b$-series described
below where $c=(r,s)$, $a=r/c$, $b=s/c$:

$a$-series:
$$
h_1^2=\pm 2bc, \ \ H\cdot h_1\equiv 0\mod bc;
$$

$b$-series:
$$
h_1^2=\pm 2ac, \ \ H\cdot h_1\equiv 0\mod ac.
$$

Then  $Y$ is isomorphic to $X$ with the isomorphism given by the
composition of the universal isomorphisms $\delta$ (if necessary),
$\nu(1,d_2)$, $T_D$ and $Tyu(\pm h_1)$ (see \eqref{isoma},
\eqref{isomb}) which have very clear geometric meaning.
\end{theorem}

Here $\delta$ is the reflection $\delta:M_X(r,H,s)\cong
M_X(s,H,r)$ (see \cite{Tyurin1}, \cite{Tyurin2}), the isomorphism
$\nu(d_1,d_2):M_X(r,H,s)\cong M_X(d_1^2r,d_1d_2H,d_2^2s)$ where
$d_1,d_2\in \bn$ and $(d_1,s)=(d_2,r)=(d_1,d_2)=1$ are some
universal isomorphisms of moduli of sheaves over an arbitrary
algebraic K3 surface. The isomorphism $T_D$ is defined by the
tensor product with $\Oc(D)$ for a class of divisors $D\in N(X)$,
thus it is also very universal. The isomorphism $Tyu(\pm
h_1):M_X(s,h_1,\pm 1)\cong X$ was geometrically defined and used
by A.N. Tyurin \cite{Tyurin2} (see also \cite{Tyurin1} and
\cite{Tyurin3}).

It follows a very interesting corollary, which shows that the
universal isomorphisms $\delta$, $\nu(d_1,d_2)$, $T_D$ and $Tyu$
are sufficient to find an isomorphism between $M_X(r,H,s)$ and
$X$, if it does exist, for a general K3 surface $X$ with $H\cdot
N(X)=\bz$. Then there exists an isomorphism between $M_X(r,H,s)$
and $X$ which is their composition.

\begin{corollary}
Let $X$ be a K3 surface with a polarization $H$ such that
$H^2=2rs$, $r,s\ge 1$, the Mukai vector $(r,H,s)$ be primitive,
and $Y=M_X(r,H,s)$ be the moduli of sheaves over $X$ with the
Mukai vector $(r,H,s)$.

Then, if $Y\cong X$ and $H\cdot N(X)=\bz$ and $X$ is general
satisfying these properties (exactly here general $X$ means that
$\rho(X)=2$ and the automorphism group of the transcendental
periods $Aut(T(X), H^{2,0}(X))=\pm 1$), there exists an
isomorphism between $Y$ and $X$ which is a composition of the
universal isomorphisms $\delta$, $\nu(d_1,d_2)$ and $T_D$ between
moduli of sheaves over $X$, and the universal Tyurin's isomorphism
$Tyu$ between a moduli of sheaves over $X$ and $X$ itself.
\label{corollarymainint}
\end{corollary}

We remark that here the isomorphisms $\delta$, $T_D$ and in
general $Tyu$ have a geometric description which does not use
Global Torelli Theorem for K3 surfaces. Only for the isomorphism
$\nu(d_1,d_2)$ we don't know a geometric construction. On the
other hand, the isomorphism $\nu(d_1,d_2)$ is very universal and
it exists even for general (with Picard number one) K3 surfaces.
By \cite{Nik0}, there exists only one isomorphism between general
algebraic K3 surfaces (or two for degree two). Thus, we can
consider this isomorphism as geometric by definition.

\section{Reminding of the main results from \cite{Mad-Nik1}, \cite{Mad-Nik2}
and \cite{Nik2}}
\label{section2}

We denote by $X$ an algebraic K3 surface
over the field $\bc$ of complex numbers. I.e. $X$ is a non-singular
projective algebraic surface over $\bc$ with the trivial canonical class
$K_X = 0$ and the vanishing  irregularity $q(X)=0$.

We denote by $N(X)$
the Picard lattice (i.e. the lattice of 2-dimensional algebraic
cycles) of $X$. By $\rho (X)=\text{rk}\ N(X)$ we denote the Picard number of
$X$. By
\begin{equation}
T(X)=N(X)_{H^2(X,\bz)}^\perp
\end{equation}
we denote the
transcendental lattice of $X$.

\par\medskip

For a Mukai vector $v=(r,c_1,s)$ where $r\in \bn$, $s \in \bz$,
and $c_1\in N(X)$, we denote by $Y=M_X(r,c_1,s)$ the moduli space
of stable (with respect to some ample $H^\prime\in N(X)$) rank $r$
sheaves on $X$ with first Chern classes $c_1$, and Euler
characteristic $r+s$.

By results of Mukai \cite{Muk1}, \cite{Muk2},
under suitable conditions on the Chern classes, the moduli space $Y$ is
always deformations equivalent
to a Hilbert scheme of 0-dimensional cycles on $X$ (of same dimension).
\par\medskip

In \cite{Mad-Nik1} for $r=s=2$, in \cite{Mad-Nik2} for $r=s=c\in
\bn$, and in \cite{Nik1} for arbitrary   $r,s\in \bn$, the
following result had been obtained.

\begin{theorem}
\label{maintheorem} Let $X$ be a K3 surface with a polarization
$H$ such that $H^2=2rs$, $r,s\ge 1$, the Mukai vector $(r,H,s)$ be
primitive, and $Y=M_X(r,H,s)$ be the moduli of sheaves on $X$ with
the Mukai vector $(r,H,s)$.

Then $Y\cong X$ if at least for one of signs $\pm$ there exists
$h_1 \in N(X)$ such that the elements $H$ and $h_1$ are contained
in a 2-dimensional sublattice $N \subset N(X)$ with $H\cdot N=\bz$
and $h_1$ belongs to either the $a$-series or the $b$-series
described below where $c=(r,s)$, $a=r/c$, $b=s/c$:

$a$-series:
\begin{equation}
h_1^2=\pm 2bc, \ \ H\cdot h_1\equiv 0\mod bc\,;
\label{aseries1}
\end{equation}

$b$-series:
\begin{equation}
h_1^2=\pm 2ac, \ \ H\cdot h_1\equiv 0\mod ac\,.
\label{bseries1}
\end{equation}

The conditions above are necessary for $H\cdot N(X)=\bz$ and
$Y\cong X$ if $\rho (X)\le 2$ and $X$ is a general K3 surface with
its Picard lattice, I. e., the automorphism group of the
transcendental periods $(T(X),H^{2,0}(X))$ of $X$ is $\pm 1$.
\end{theorem}

{\bf Remark 2.1.} In \cite{Nik1}, in \eqref{aseries1},
\eqref{bseries1} some additional primitivity conditions on $h_1$
were required. Thus, the corresponding conditions in \cite{Nik1} were

for $a$-series:
\begin{equation} h_1^2=\pm 2bc, \ \ H\cdot h_1\equiv
0\mod bc,\ \ H\cdot h_1\not\equiv 0\mod bcl_1,\ \ h_1/l_2\notin
N(X),
\label{aseries2}
\end{equation}
where $l_1$ and $l_2$ are any primes such that $l_1^2|a$ and
$l_2^2|b$;

for $b$-series:
\begin{equation}
h_1^2=\pm 2ac, \ \ H\cdot h_1\equiv 0\mod ac,\ \ H\cdot
h_1\not\equiv 0\mod acl_1,\ \ h_1/l_2\notin N(X),
\label{bseries2}
\end{equation}
where $l_1$ and $l_2$ are any primes such that $l_1^2|b$ and
$l_2^2|a$.

\medskip

Let us show that actually \eqref{aseries1} is equivalent to
\eqref{aseries2}, and \eqref{bseries1} is equivalent to
\eqref{bseries2} (where we assume that $H\cdot N(X)=\bz$). 
Thus, in sufficient part of Theorem 
\ref{maintheorem} it is better to use \eqref{aseries1},
\eqref{bseries1}, and in necessary part of Theorem
\ref{maintheorem} it is better to use \eqref{aseries2},
\eqref{bseries2}.

According to \cite{Nik1}, for some $\mu\mod{2abc^2}\in
(\bz/2abc^2)^\ast$ and $d\in \bn$ such that $\mu^2\equiv d\mod
4abc^2$, \eqref{aseries1} is equivalent to the system of
diophantine conditions (see (3.54) in \cite{Nik1})
\begin{equation}
p^2-dq^2=\pm 4ac,\ \ p\equiv \mu q\mod 2ac,
\label{aseries3}
\end{equation}
and \eqref{aseries2} is equivalent to
\begin{equation}
p^2-dq^2=\pm 4ac,\ \ p\equiv \mu q\mod 2ac, \ \ (a,p,q)=1,\ \
p\not\equiv \mu q\mod 2asl  \label{aseries4}
\end{equation}
for any prime $l|b$.

Let us show that \eqref{aseries3} implies \eqref{aseries4}. We can
find $\mu_0\mod 2a^2c^2$ such that $\mu_0\equiv \mu\mod 2ac^2$ and
$d\equiv \mu_0^2\mod 4a^2c^2$. We then get $p^2-dq^2\equiv
(p-\mu_0 q)(p+\mu_0 q)\equiv \pm 4ac\mod 4a^2c^2$. Then
$$
\frac{p-\mu_0q}{2ac}\cdot (p+\mu_0 q)\equiv \pm 2\mod 2ac
$$
where $(p-\mu_0 q)/(2ac)$ is an integer. This leads to a
contradiction if a prime $l|(a,p,q)$.

Similarly, since $\mu^2\equiv d\mod 4abc^2$, we obtain that
$$
\frac{p-\mu q}{2ac}\cdot \frac{p+\mu q}{2}\equiv \pm 1\mod bc.
$$
This leads to a contradiction if $p-\mu q\equiv 0\mod 2acl$ for a
prime $l|b$.

Thus, \eqref{aseries3} implies \eqref{aseries4}, and they are
equivalent indeed.

\medskip

The sufficient part of the proof of Theorem \ref{maintheorem} in
\cite{Nik1} and similar Theorems in \cite{Mad-Nik1},
\cite{Mad-Nik2} used global Torelli Theorem for K3 surfaces
\cite{PShShaf}. I. e., under conditions of Theorem
\ref{maintheorem}, we proved that the K3 surfaces $X$ and $Y$ have
isomorphic periods. By global Torelli Theorem,  then $X$ and $Y$
are isomorphic.

In Sect. \ref{section3} below we will give a geometric
construction of the isomorphism between $X$ and $Y$ which is
similar to our considerations in \cite{Mad-Nik3}. This is the main
result of this paper.

\section{Geometric interpretation of the main results from 
\cite{Mad-Nik1}, \cite{Mad-Nik2} and
\cite{Nik2}} \label{section3}

We use notations of the previous Section \ref{section2}.

We shall use the natural isomorphisms between moduli spaces of sheaves over a K3
surface $X$.

\begin{lemma}
\label{lemmadelta} Let $(r,H,s)$ be a primitive Mukai vector for a
K3 surface $X$ and $r,s\ge 1$. Then one has an isomorphism, called
{\bf reflection}, see \cite{Tyurin2},
$$
\delta:M_X(r,H,s)\cong M_X(s,H,r).
$$
 \end{lemma}

In (\cite{Tyurin1}, (4.11)) and \cite{Tyurin2}, Lemma 3.4), the geometric construction of the reflection 
$\delta$ is given under the condition that $M_X(r,H,s)$ contains
a regular bundle. On the other hand, by global Torelli
Theorem for K3 surfaces \cite{PShShaf}, existence of such
isomorphism is obvious. See similar proof of Theorem
\ref{maintheorem1} below. 

\begin{lemma}
\label{lemmaTD}
Let $(r,H,s)$ be a Mukai vector for a K3 surface $X$ and $D\in N(X)$.
Then one has the natural isomorphism given by the tensor product
$$
T_D:M_X(r,H,s)\cong M_X(r,H+rD,s+r(D^2/2)+D\cdot H),\ \  \E\mapsto \E\otimes \Oc(D).
$$
Moreover, here Mukai vectors
$$
v=(r,H,s),\ \ v_1=(r,H+rD,s+r(D^2/2)+D\cdot H)
$$
have the same general common divisor and the same square under Mukai pairing.
In particular, they are primitive and isotropic simultaneously.
 \end{lemma}

\medskip

We also use the isomorphisms between moduli of sheaves over a K3
surface and the K3 surface itself found by A.N.Tyurin in 
(\cite{Tyurin2}, Lemma 3.3).

\begin{lemma}
\label{lemmaTyu}
For a K3 surface $X$ and $h_1\in N(X)$ such that $\pm h_1^2>0$
and
\begin{equation}
h^0\Oc(h_1)=h^0\Oc(-h_1)=0\  if\  h_1^2<0,
\label{Tyurin1}
\end{equation}
there is a geometric Tyurin's isomorphism
$$
Tyu(\pm h_1):M_X(\pm h_1^2/2,h_1, \pm 1)\cong X .
$$
\end{lemma}

Like in Theorem \ref{maintheorem1} below, using global Torelli
Theorem for K3 surfaces \cite{PShShaf}, one can show that even
without Tyurin's condition \eqref{Tyurin1} always there exists
some isomorphism
$$
M_X(\pm h_1^2/2,h_1, \pm 1)\cong X.
$$
We call such isomorphisms also as {\it Tyurin's isomorphisms.}
When $h_1$ also satisfies \eqref{Tyurin1}, there exists a direct
Tyurin's geometric construction of some isomorphism $M_X(\pm
h_1^2/2,h_1,\pm 1)\cong X$.

\medskip

We shall use the following result which had been proved implicitly
in \cite{Nik2}. See Sect. 2.3 and the proof of Theorem 2.3.3 in
\cite{Nik2}. In \cite{Nik2} much more difficult results related to
arbitrary Picard lattice had been considered, and respectively the
proofs were long and difficult. Therefore, below we also give a
much simpler proof of this result.

\begin{theorem}
Let $v=(r,H,s)$  be an isotropic Mukai vector on a K3 surface $X$
where $r,\ s \in \bn$, $H\in N(X)$, $H^2=2rs$,  and $H$ is
primitive (then $v$ is also primitive).

Let $d_1,d_2\in \bn$ and
$(d_1,s)=(d_2,r)=(d_1,d_2)=1$.

Then the Mukai vector
$v_1=(d_1^2r, d_1d_2H,d_2^2s)$
is also primitive, and  there exists a natural isomorphism of moduli of sheaves
$$
\nu(d_1,d_2): M_X(r,H,s)\cong M_X(d_1^2r,d_1d_2H,d_2^2s).
$$
\label{maintheorem1}
\end{theorem}

\begin{proof} Let us consider the case $(d_1,d_2)=(d,1)$ where $(d,s)=1$ (general case
is similar).

Let $c=(r,s)$ and $a=r/c$, $b=s/c$. Then $v=(ac,H,bc)$ where
$H^2=2abc^2$ and $H$ is primitive. By global Torelli Theorem
\cite{PShShaf}, it is enough to show that periods of
$Y=M_X(r,H,s)$ and $Y_1=M_X(d^2r,dH,s)$ are isomorphic.

By results of Mukai \cite{Muk1}, \cite{Muk2}, cohomology of $Y$
(and similarly of $Y_1$) are equal to
$$
H^2(Y,\bz)=v^\perp/\bz v
$$
where we consider $v$ as the element of Mukai lattice $H^\ast(X,\bz)$, and
$H^{2,0}(Y)$ is the image of $H^{2,0}(X)$.

Using the variant of Witt's Theorem from \cite{PShShaf}, we can model the necessary
calculations as follows.

Let $U^{(1)}$ be an even unimodular hyperbolic plane  with the
basis $e_1,\ e_2$ where $e_1^2=e_2^2=0$ and $e_1\cdot e_2=-1$. Let
$U^{(2)}$ be another even unimodular hyperbolic plane with the
bases $f_1,\ f_2$ where $f_1^2=f_2^2=0$ and $f_1\cdot f_2=1$. We
consider the orthogonal sum $U^{(1)}\oplus U^{(2)}$ (the model of
$H^\ast (X,\bz)$; to get $H^\ast(X,\bz)$, one should add to
$U^{(1)}\oplus U^{(2)}$ an unimodular even lattice of signature
$(2, 18)$ which is the same for all three $X,\,Y,\,Y_1$).  Then
$$
v=ace_1+bce_2+H,\ \ H=abc^2f_1+f_2,
$$
thus
$$
v=ace_1+bce_2+abc^2f_1+f_2.
$$
Then $N(X)=\bz H$ models the Picard lattice of $X$ and $T(X)=\bz t$, $t=-abc^2f_1+f_2$
models the transcendental lattice of $X$.

We have $\xi=x\,e_1+y\,e_2+z\,f_1+w\,f_2\in v^\perp$ if and only
if $-bcx-acy+z+abc^2w=0$, equivalently $z=bcx+acy-abc^2w$ where
$x,y,w,z\in \bz$. Thus,
$$
\xi=x(e_1+bc\,f_1)+y(e_2+ac\,f_1)+w(-abc^2\,f_1+f_2)\,
$$
and $\alpha=e_1+bc\,f_1$, $\beta=e_2+ac\,f_1$, $t=-abc^2\,f_1+f_2$
give the basis of $v^\perp$. We have $v=ac\,\alpha+bc\,\beta+t$
and
$$
ac\,\alpha\hskip-7pt\mod {\bz v}+bc\,\beta\hskip-7pt\mod{\bz v}+
t\hskip-7pt\mod{\bz v}=0.
$$
It follows that $\overline{\alpha}=\alpha\hskip-5pt\mod {\bz v}$,
$\overline{\beta}= \beta\hskip-5pt\mod{\bz v}$ give a basis of
$v^\perp/\bz v$ which models $H^2(Y,\bz)$. We have
$\overline{\alpha}^2=\overline{\beta}^2=0$ and
$\overline{\alpha}\cdot \overline{\beta}=-1$. Thus $v^\perp/\bz
v\cong U$. We see that $\overline{t}=t\hskip-5pt\mod{\bz
v}=-ac\overline{\alpha}-bc\overline{\beta}$ and then
$\widetilde{t}=\overline{t}/c=-a\overline{\alpha}-b\overline{\beta}\in
(v^\perp/\bz v)$, and $\bz \widetilde{t}$ models the
transcendental lattice $T(Y)$ of $Y$. Its orthogonal complement
$\bz h$ where $h=a\overline{\alpha}-b\overline{\beta}$ models the
Picard lattice $N(Y)$ of $Y$. We have $h^2=2ab >0$.

Let us make similar calculations for $Y_1$. We have
$$
v_1=d^2ac\,e_1+bc\,e_2+dH,\ \ H=abc^2\,f_1+f_2,
$$
and
$$
v_1=d^2ac\,e_1+bc\,e_2+dabc^2\,f_1+d\,f_2.
$$
We have $\xi_1=x_1e_1+y_1e_2+z_1f_1+w_1f_2\in v_1^\perp$ if and
only if $-bcx_1-d^2acy_1+dz_1+dabc^2w_1=0$. Since $(d,bc)=1$, we
obtain $x_1=d\widetilde{x}_1$, and
$z_1=bc\widetilde{x}_1+dacy_1-abc^2w_1$ where
$\widetilde{x}_1,y_1,w_1,z_1\in \bz$. Then
$$
\xi_1=\widetilde{x}_1(de_1+bcf_1)+y_1(e_2+dacf_1)+w_1(-abc^2f_1+f_2)
$$
and $\alpha_1=d\,e_1+bc\,f_1$, $\beta_1=e_2+dac\,f_1$,
$t=-abc^2\,f_1+f_2$ give the basis of $v_1^\perp$. We have
$v_1=dac\,\alpha_1+bc\,\beta_1+d\,t$, and
$$
dac\,\alpha_1\hskip-7pt\mod {\bz v_1}+bc\,\beta_1\hskip-7pt\mod{\bz v_1}+
d\,t\hskip-7pt\mod{\bz v_1}=0.
$$
Since $(d,bc)=1$, we see that $t\hskip-5pt\mod{\bz v_1}=c\tilde{t}_1$,
$\beta_1\hskip-5pt\mod{\bz v_1}=d\tilde{\beta}_1$ where
$\tilde{t}_1, \tilde{\beta}_1\in v_1^\perp/\bz v_1$, and
$$
a\,\alpha_1\hskip-7pt\mod {\bz v_1}+b\tilde{\beta}_1+\tilde{t}_1=0.
$$
We have  $(\alpha_1\hskip-5pt\mod {\bz v_1})^2=\tilde{\beta}_1^2=0$ and
$(\alpha_1\hskip-5pt\mod {\bz v_1})\cdot \tilde{\beta}_1=-1$.
Thus, again $\alpha_1\hskip-5pt\mod {\bz v_1}$, $\tilde{\beta}_1$ give canonical
generators of the unimodular lattice $U$.
Then they give a basis of $v_1^\perp/\bz v_1$ which models $H^2(Y_1)$. Moreover
$\bz \tilde{t}_1$ where
$\tilde{t}_1= -a\,\alpha_1\hskip-7pt\mod {\bz v}-b\,\tilde{\beta}_1$ (from above),
 models the transcendental lattice of $Y_1$.
Its orthogonal complement $\bz h_1$ where
$h_1=a\,\alpha_1\hskip-5pt\mod {\bz v}-b\,\tilde{\beta}_1$ models the
Picard lattice of $Y_1$.

We see that our descriptions for $Y$ and $Y_1$ above are evidently
isomorphic if we identify $\overline{\alpha}$, $\overline{\beta}$
with $\alpha_1\hskip-5pt\mod {\bz v}$ and $\tilde{\beta}_1$
respectively. This shows that $Y$ and $Y_1$ have isomorphic
periods and are isomorphic by global Torelli Theorem for K3
surfaces \cite{PShShaf}.

This finishes the proof.
\end{proof}

\medskip

{\bf Remark 3.1.} It would be very interesting to find a direct
geometric proof of Theorem \ref{maintheorem1} which does not use
global Torelli Theorem for K3 surfaces. It seems, considerations
by Mukai in  \cite{Muk3} are related with this problem 
(especially see Theorem 1.2, in  \cite{Muk3}). 

On the other hand, the isomorphism $\nu(d_1,d_2)$ is very universal, and
it exists even for general (with Picard number one) K3 surfaces.
By \cite{Nik0}, there exists only one isomorphism (or two isomorphisms 
for the degree two) between algebraic K3 surfaces with Picard number one. 
Thus, we can consider this isomorphism as geometric by definition.

\medskip

We will show below that under the conditions of Theorem
\ref{maintheorem}, there exists an isomorphism between $X$ and $Y$
which is the composition of the universal geometric isomorphisms above.

\begin{theorem}
\label{maintheorem2} Let $X$ be a K3 surface with a polarization
$H$ such that $H^2=2rs$, $r,s\ge 1$, the Mukai vector $(r,H,s)$ be
primitive, and $Y=M_X(r,H,s)$ be the moduli of sheaves on $X$ with
the Mukai vector $(r,H,s)$.

Assume that at least for one of signs $\pm$ there exists $h_1 \in
N(X)$ such that the elements $H$ and $h_1$ are contained in a
2-dimensional sublattice $N \subset N(X)$ with $H\cdot N=\bz$ and
$h_1$ belongs to either the $a$-series or the $b$-series described
below where $c=(r,s)$, $a=r/c$, $b=s/c$:

$a$-series:
$$
h_1^2=\pm 2bc, \ \ H\cdot h_1\equiv 0\mod bc;
$$

$b$-series:
$$
h_1^2=\pm 2ac, \ \ H\cdot h_1\equiv 0\mod ac.
$$

Then  $Y$ is isomorphic to $X$ with the isomorphism given by the
composition of the reflection $\delta$ (if $h_1$ belongs to the
$a$-series), $\nu(1,d_2)$ (for some $d_2$), $T_D$ (for some $D\in
N$) and $Tyu(\pm h_1)$ (see \eqref{isoma}, \eqref{isomb}) where
$d_2$ and $D$ are defined in the proof below.
\end{theorem}

\begin{proof} We use the description of the pair $H\in N$ given in Proposition
3.1 in \cite{Nik2} which directly follows from $\text{rk\ } N=2$
and $H\cdot N=\bz$. We denote $d=-\det N$ and $\bz\delta$ the
orthogonal complement to $H$ in $N$. Then $\delta^2=-2abc^2d$.

We have
\begin{equation}
N=[H,\delta, w=\frac{\mu H+\delta}{2abc^2}],\
\mu\mod{2abc^2} \in (\bz/2abc^2)^\ast,\ \mu^2\equiv d\mod 4abc^2
\label{latticeN1}
\end{equation}
where $[\ \ \cdot\  \ ]$ means `generated by $\cdot$ ". Here $d$
and $\pm \mu\mod2{abc^2}\subset (\bz/2abc^2)^\ast$ give the
complete invariants of the pair $H\in N$ up to isomorphisms.

It follows that
\begin{equation}
N=\{z=\frac{xH+y\delta}{2abc^2}\ |\ x,y\in \bz\ \text{and}\ x\equiv \mu y\mod 2abc^2\}
\label{latticeN2}
\end{equation}
where
\begin{equation}
z^2=\frac{x^2-dy^2}{2abc^2}\ .
\label{latticeN3}
\end{equation}

In the conditions of Theorem \ref{maintheorem1}, let us assume that $h_1$ belongs to
the $a$-series.

Then $h_1=(rH+s\delta)/(2abc^2)$ where $r,s\in \bz$ and $r\equiv
\mu s\mod{2abc^2}$. Since $h_1^2=\pm 2bc$, it follows that
$r^2-ds^2=\pm (2bc)(2abc^2)$. We have $H\cdot h_1=r\equiv
0\mod{bc}$. Since $r\equiv \mu s\mod {2abc^2}$, it follows
$s\equiv 0\mod{bc}$. Denoting $r=pbc$ and $s=qbc$ where $p,q\in
\bz$, we get that $h_1=(pH+q\delta)/(2ac)$ where $p\equiv \mu
q\mod{2ac}$ and $p^2-dq^2=\pm 4ac$.

We have $\delta=2abc^2\,w-\mu H$. Then
$$
h_1=\frac{pH+q(2abc^2 w-\mu H)}{2ac}=\frac{p-\mu q}{2ac} H+qbc w\equiv
\frac{p-\mu q}{2ac} H\mod{bcN}.
$$

Since $p\equiv \mu q\mod{2ac}$, there exists $d_2\in \bn$ such that 
\begin{equation}
d_2\equiv
\frac{p-\mu q}{2ac}\mod bc\,. 
\label{formula}
\end{equation} 
Then
$$
h_1=d_2H+bc D, \ \ D\in N.
$$
Then (see Lemma \ref{lemmadelta})
$$
\delta:(ac,H,bc)\to (bc,H,ac),
$$
and (see Lemma \ref{lemmaTD})
$$
T_D:(bc,d_2H,d_2^2ac )\to (bc, h_1,\pm 1)
$$
since $(bc,d_2H,d_2^2ac)$ is isotropic Mukai vector and $h_1^2=\pm
2bc$. Since $(bc,h_1,\pm 1)$ is evidently primitive,
$(bc,d_2H,d_2^2ac)$ is also primitive, and then $(d_2,bc)=1$. By
Theorem \ref{maintheorem1}, then
$$
\nu (1,d_2):M_X(bc,H,ac)\cong M_X(bc,d_2H,d_2^2ac)
$$
and
$$
T_D:M_X(bc,d_2H,d_2^2ac)\cong M_X(bc,h_1,\pm 1).
$$
At last (see Lemma \ref{lemmaTyu}),
$$
Tyu(\pm h_1):M_X(bc,h_1,\pm 1)\cong X.
$$
Thus, we obtain the natural isomorphism
\begin{equation}
Tyu(\pm h_1)\cdot T_D\cdot \nu(1,d_2)\cdot \delta: M_X(r,H,s)\cong
X . \label{isoma}
\end{equation}

If $h_1$ belongs to the $b$-series, similarly we can show that
$h_1=d_2H+acD$ for $d_2\in \bn$ and $D\in N$, and we obtain the natural isomorphism
\begin{equation}
Tyu(\pm h_1)\cdot T_D\cdot \nu(1,d_2): M_X(r,H,s)\cong X .
\label{isomb}
\end{equation}
We don't need the reflection $\delta$ in this case.

This finishes the proof of Theorem \ref{maintheorem2}.
\end{proof}

Since conditions of Theorem \ref{maintheorem2} are also necessary
for a general K3 surface $X$ with $H\cdot N(X)=\bz$ and 
$Y\cong X$ (see Theorem \ref{maintheorem}), we
also obtain an interesting

\begin{corollary}
Let $X$ be a K3 surface with a polarization $H$ such that
$H^2=2rs$, $r,s\ge 1$, the Mukai vector $(r,H,s)$ be primitive,
and $Y=M_X(r,H,s)$ be the moduli of sheaves over $X$ with the
Mukai vector $(r,H,s)$.

Then, if $Y\cong X$ and $H\cdot N(X)=\bz$ and $X$ is general
satisfying these properties (exactly here general $X$ means that
$\rho(X)=2$ and the automorphism group of the transcendental
periods $Aut(T(X), H^{2,0}(X))=\pm 1$), there exists an
isomorphism between $Y$ and $X$ which is a composition of the
universal isomorphisms $\delta$, $\nu(d_1,d_2)$ and $T_D$ between
moduli of sheaves over $X$, and the universal Tyurin's isomorphism
$Tyu$ between the moduli of sheaves over $X$ and $X$ itself.
\label{corollarymain}
\end{corollary}

Here it is important that the isomorphisms $\delta$, $T_D$ and 
$Tyu$ have a geometric description which does not use
Global Torelli Theorem for K3 surfaces. Only for the isomorphism
$\nu(d_1,d_2)$ we don't know a geometric construction. On the
other hand, the isomorphism $\nu(d_1,d_2)$ is very universal and is 
geometric by definition. See our considerations at the beginning of 
the Section.

\

\

C.G.Madonna \par
Math. Dept.,
CSIC, C/ Serrano 121,
28006 Madrid,
SPAIN

carlo@madonna.rm.it \ \
cgm@imaff.cfmac.csic.es

\

\

V.V. Nikulin \par
Deptm. of Pure Mathem. The University of Liverpool, Liverpool\par
L69 3BX, UK;
\vskip1pt
Steklov Mathematical Institute,\par
ul. Gubkina 8, Moscow 117966, GSP-1, Russia

vnikulin@liv.ac.uk \ \
vvnikulin@list.ru

\end{document}